\documentclass[a4paper, twoside,10pt,fullpage]{article}
\usepackage{fontenc}
\usepackage{amsmath}
\usepackage{amssymb}
\usepackage{amsfonts}
\usepackage{amsthm}
\usepackage{graphicx}
\usepackage{eucal}
\usepackage{mathrsfs}
\usepackage{amsthm}
\usepackage{verbatim}
\usepackage{epsfig}
\usepackage{hyperref}
\usepackage{float}
\usepackage{color}

\newtheorem{thm}{Theorem}

\newtheorem{cor}[thm]{Corollary}

\newtheorem{que}{Question}
\newtheorem{prob}[que]{Problem}

\newcommand{\R}{{\mathbb{R}}}

\newcommand{\Z}{{\mathbb{Z}}}

\def\H{{\cal{H}}}
\def\ep
{\hfill{$\Box$}}

\def\lap{\bigtriangleup}

\def\gam{\gamma}
\def\EL{\rm EL}
\def\length{\mbox{length}}
\def\Gam{\Gamma}
\def\supp{\mbox{supp}}
\def\diam{\mbox{diam}}
\def\MD{\Omega}

\begin{document}
\title{Lack of Sphere Packing of Graphs via Non-Linear Potential Theory}
\author{Itai Benjamini \and Oded Schramm  \footnote{Oded died while solo climbing Guye Peak in Washington State on September 1, 2008.}}
\maketitle

\begin{abstract}
It is shown that there is no quasi-sphere packing of
the lattice grid
$\Z^{d+1}$ or a co-compact hyperbolic lattice of $\H^{d+1}$ or the $3$-regular tree $ \times  \mbox{  } \Z$, in $\R^d$, for all $d$.
A similar result is proved for some other graphs too.
Rather than %
using a direct geometrical approach,
the main tools we are using are from non-linear potential theory.

\end{abstract}

\section{Introduction}
In this paper we would like to study sphere (and quasi-sphere) packing of
graphs in Euclidean spaces.
Suppose that $P=(P_v:v\in V)$ is an indexed packing of spheres.  This just
means that $V$ is some set, to each $v$ corresponds
a sphere $P_v\subset\R^n$, and the interiors of the spheres are disjoint.
Let $G$ be the graph with vertices $V$ such that there is an edge joining
$v$ and $u$ iff $P_v$ and $P_u$ are tangent.  Then $G$ is the
{\it contacts graph\/} of $P$.  (There are no multiple edges in $G$).
In a quasi-sphere packing, We allow $P_v$ to be any domain for which,
the ratio between the outer radius and the inner radius, are uniformly
bounded, for all $v \in V$.
The two dimensional case is well understood, if $G$ admits a circle
packing, then
It is easy to see that $G$ is planar, the circle packing theorem provides
a converse:

\begin{thm} (Circle Packing Theorem)
\label{cpt}
Let $G=(V,E)$ be a finite planar graph with no
loops or multiple edges, then there is a disk packing
$P=(P_v:v\in V)$ in $\R^2$ with contacts graph $G$.
\end{thm}

This theorem was first proved by Koebe~\cite{ko}.
Recently, at least $7$ other proofs have been found; some
of the more accessible ones can be found
in~\cite{m-r}, ~\cite{co}, ~\cite{b-s}.
\medskip

Which graphs admit a quasi-sphere packing in $\R^d, d>2$ ?
We know of two kinds of invariants which imply that the graph cannot be
packed in $\R^d$. One is the separation function of the graph. This invariant
is geometrical, and in particular it can be used to show that, the lattice
$\Z^{d+1}$ can not be quasi-sphere packed in $\R^d$. However this invariant is
not good enough to show, for instance, that a co-compact lattice in the $d+1$ dimensional hyperbolic space $\H^{d+1}$,
or the $3$-regular tree $ \times \mbox{  }  \Z$ cannot be packed in $\R^d$.
Note that the $3$-regular tree $ \times \mbox{  }  \Z$ admits square root separators, similar to planar graphs,
although is not planar,
see~\cite{bs:har}.
We will not discuss separation here, regarding
many properties of the separation function
for graphs, including the application above, see~\cite{bs:sep}.

In this note we will use another invariant, which is  formulated in the
language of non-linear potential theory.  In particular it gives,

\begin{thm}
\label{packZH}
The lattice $\Z^{d+1}$ as well as any co-compact lattice in $\H^{d+1}$
or $3$-regular tree $ \times \mbox{  } \Z$
can not be quasi-sphere packed in $\R^d$.
\end{thm}

\noindent
This theorem is a corollary of our main result, which is more technical.

\begin{thm}
\label{packgen}
Assume $G$ admits a uniform bound on its %
vertex degree. If $G$ is
$d$-nonparabolic, yet admits no non constant, $d$-Dirichlet, $d$-harmonic functions,
then $G$ has no quasi-sphere packing in $\R^d$.
\end{thm}
\medskip

\noindent
All the notions in theorem~\ref{packgen} are defined in the next section.
The rest of the sections contains the proof of this theorem, concluding
with the derivation of Theorem~\ref{packZH}. We end with a couple of problems.
\medskip

\noindent
{\bf Remarks:}
\medskip

\noindent
1) Existence of quasi-sphere packing is a rough-isometry invariant of the
graph.
\medskip

\noindent
2) Kuperberg and Schramm~\cite{ku-sc}
gave bounds on the possible average kissing numbers of sphere packing in
$\R^3$. See~\cite{co-sl} for more on fixed radius sphere packing.
\medskip

\noindent
3)  $\H^2 \times \R$ is a bounded geometry %
Riemannian metric on $\R^3$. For the reasons
$3$-regular tree $ \times \mbox{ } \Z$ is not sphere packable, it is possible to show that there is no packable graph
in Euclidean $\R^d$ roughly isometric to $\H^2 \times \R$. Which bounded geometry %
Riemannian metrics on $\R^3$
admit a roughly %
isometric packable graph?

\section{Notations and Terminology}

Let $G=(V,E)$ be a graph.
For convenience, we usually only consider graphs
with no loops or multiple edges (but the results do
apply to multigraphs).
We sometimes use $\{v,u\}$ to denote the edge with endpoints $v,u$.

The graphs we shall consider will be connected
and locally finite.  The latter means that the number of edges
incident with any particular vertex is finite.
\medskip

We now bring some useful definitions from potential theory.
\medskip

\noindent
{\bf Definition.} The {\it $p-$Dirichlet energy\/} of a function $f:V\to\R$
is defined by
$$
D_p(f)= 1/2 \sum_{\{v,u|\{u,v\} \in E\}} |f(u)-f(v)|^p
$$

A {\it $p$-Dirichlet function\/} is an $f:V\to\R$ with $D_p(f)<\infty$.
The space of all $p$-Dirichlet functions on $G$ is denoted by $D_p(G)$.
\medskip

\noindent
{\bf Definition.} A graph $G$ is {\it $p$-parabolic} if
$$
\inf \{ D_p(f) | f\mbox{ has finite support, } f(v_0)=1\} =0 \mbox{ for some
vertex } v_0
$$

Otherwise we say that the graph is {\it $p$-nonparabolic}.
Nonparabolic is called sometimes hyperbolic. $2$-parabolic has a
probabilistic meaning, simple random walk on the graph is recurrent.
Soardi and Yamasaki~\cite{so-ya} proved assuming a uniform bound on the degrees that the parabolic index
($ = \sup\{p \geq 1| G \mbox{ is p-parabolic }\}$), is invariant under rough
isometries, (see also~\cite{so:1} and ~\cite{ho-so:2}). Maeda~\cite{ma}
proved that the parabolic index of $\Z^d$ is exactly $d$,
(see also~\cite{so-ya} and~\cite{ho:1}). It is not hard to
see that the binary tree is $p$-nonparabolic for any $p \geq 1$.

For $f:V\to\R$, and $1 < p < \infty$, set
$$
\lap_p f(v)= \sum_{\{u| \{u,v\} \in E\}} |f(u)-f(v)|^{p-2}(f(u)-f(v)). %
$$

\noindent
{\bf Definition.} A function  $f:V\to\R$ is called {\it $p$-harmonic} if
$\lap_p f(v)= 0$ at every $v \in V$.
\medskip

A {\it metric\/} $m$ on a graph $G=(V,E)$ is a positive function
$m:E\to (0,\infty)$.

\noindent
The {\it gradient\/} of a function $f:V\to \R$
with respect to a metric $m$ is defined by
$$
|df(e)/m(e)| =  |(f(v)-f(u))/m(e)|,
$$
\noindent
where $\{v,u\}=e$.

The {\it natural metric\/} on $G$ is the metric where each
edge gets %
weight $1$.
In the absence of another metric,
all metric related notions are assumed to be with respect to
the natural metric.

Two metrics $m,m'$ are {\it mutually bilipschitz\/},
if the ratios $m/m'$ and $m'/m$ are bounded.

Let $G=(V,E)$ be a connected, locally finite graph,
and let $m$ be a metric on $G$.
The $m$-length of a path $\gam$ in $G$ is the sum of
$m(e)$ over all edges in $\gam$,
$$
\length_m(\gam)=\sum_{e\in\gam} m(e).
$$
We define the $m$-distance $d_m(v,u)$ between any two vertices
$v,u\in V$ to be the infimum of the $m$-lengths of paths connecting
$v$ and $u$.  Then $(V,d_m)$ is a metric space.

\section{$p$-resolving metrics}

Below, we introduce the notion of a $p$-resolvable graph, and will
see that a $p$-nonparabolic $p$-resolvable graph has non-constant,
$p$-harmonic functions.
\medskip

\noindent
Let $G=(V,E)$ be some graph, and
let $\Gam$ be a collection of (infinite) paths in $G$.
Then $\Gam$ is {\it $p$-null\/} if there is an $L^p(E)$ metric
on $G$ such that $\length_m(\gam)=\infty$ for every $\gam\in\Gam$.

It is easy to see that $\Gam$ is $p$-null iff its {\it $p$-extremal length\/}
$$
\EL(\Gam)=\sup_m\inf_{\gam\in\Gam} \frac{\length_m(\gam)^p}{\|m\|^p},
$$
is infinite.
(Extremal length was imported to the discrete setting by
Duffin~\cite{du}. See~\cite{so:2} for more about extremal
length on graphs.)
When $\Gam$ is a collection of paths and a property holds for every
$\gam\in\Gam$, except for a $p$-null family, we shall say that the
property holds for {\it $p$-almost every\/} $\gam\in\Gam$.

Let $m$ be a metric on $G$, and recall that $d_m$ is
the associated distance function.
Let $C_m(G)$ denote the completion of $(V,d_m)$,
and let the $m$-{\it boundary} of $G$ be $\partial_mG =C_m(G)-V$.
We use $d_m$ to also denote the metric of the completion $C_m(G)$.

The metric $m$ will be called {\it $p$-resolving\/} if
it is in $L^p(E)$ and for every $x\in \partial_mG$ the
collection of half infinite paths in $G$ that converge to
$x$ in $C_m(G)$ is $p$-null.
$G$ is {\it $p$-resolvable\/} if it has a $p$-resolving metric.

Note that if $m$ is a $p$-resolving metric and $m'$ is another $L^p$
metric satisfying $m'\geq m$, then $m'$ is also $p$-resolving.

Theorem~\ref{resolving} below
shows that any $p$-parabolic graph is $p$-resolvable,
in fact, there is a metric $m$ with $\partial_m G=\emptyset$.
On the other hand, the next theorem shows that a $p$-nonparabolic
graph with no non-constant, $p$-harmonic functions, is not $p$-resolvable,
for example,
$\Z^3$ or any co-compact lattice in hyperbolic $d$-space $d>2$ is not
$2$-resolvable.
Any bounded valence planar graph is $2$-resolvable ~\cite{bs:har}.

\begin{thm}
\label{resolving}
Let $G=(V,E)$ be an infinite, connected, locally finite, $p$-resolvable graph.
If $G$ is $p$-nonparabolic, then there are non-constant, $p$-Dirichlet %
$p$-harmonic functions on $G$.
\end{thm}

We shall need the following results.

\begin{thm} (Yamasaki~\cite{ya:1})
\label{recel}
Let $G$ be a locally finite connected graph, and let $\Gam$ be
the collection of all infinite paths in $G$.
Then $G$ is $p$-parabolic if and only if $\Gam$ is $p$-null.
\end{thm}

(Yamasaki~\cite{ya:1} considers only paths that start at a fixed base vertex,
but this is equivalent.)


\begin{thm} (Yamasaki~\cite{ya:2})
\label{hde}
There are non-constant, $p$-Dirichlet, %
$p$-harmonic, functions
on $G$ if and only if there is an $f\in D^p(G)$ such that
for every $c\in\R$ the collection of all one-sided-infinite paths $\gam$
in $G$ with $\lim_n f(\gam(n))\neq c$ is not $p$-null.
\end{thm}

\noindent
{\bf Proof of theorem~\ref{resolving}:}
Assume that $G$ is $p$-nonparabolic, and $m$ is a
$p$-resolving metric on $G$.
Let $\Gam$ be the collection of all infinite paths
$\gam=(\gam(0),\gam(1),\dots)$ in $G$.
Almost all paths $\gam$ in $\Gam$ have a limit $\lim_n\gam(n)$
in $\partial_mG$, since the $m$-length of those that do not is infinite.
(The limit is in the metric $d_m$, of course.)

We now define $ \supp(\Gam)$, the {\it support\/} of $\Gam$ in $\partial_mG$,
as the intersection of all closed sets $Q\subset \partial_m G$ such
that for almost every $\gam\in \Gam$ the limit $\lim_n\gam(n)$ is
in $Q$.  Because there is a countable basis for the topology
of $\partial_m G$, and a countable union of $p$-null collections of curves is
$p$-null,
almost every $\gam\in \Gam$ satisfies $\lim_n\gam(n)\in \supp(\Gam)$.

Since $G$ is $p$-nonparabolic, we know from theorem~\ref{recel} that
the extremal length of $\Gam$ is finite.
Consequently, $\supp(\Gam)$ is not empty.
Moreover, the assumption that $m$ is $p$-resolving
shows that $\supp(\Gam)$ consists of more than a single point.
Let $x_0$ be an arbitrary point in $\supp(\Gam)$.
Define $f:V\to\R$ by setting $f(p)=d_m(x_0,p)$.
It is clear that $|df(e)|\leq m(e)$ holds for $e\in E$.
Consequently, $f\in D^p(G)$.

Pick some $\delta>0$ that is smaller than the $d_m$-diameter
of $\supp(\Gam)$.  Consider the set
$A_\delta=\{x\in\supp (\Gam): d(x_0,x)<\delta\}$,
and let $\Gam_\delta$ be the set of $\gam\in \Gam$ such
that $\lim_n\gam(n)\in A_\delta$.
Since $\supp(\Gam)$ is not contained in $\bar{A}_\delta$
or in $\partial_m G-A_\delta$,
{}from the definition of $\supp(\Gam)$ it follows that both
$\Gam_\delta$ and $\Gam-\Gam_\delta$ are not $p$-null.
For every $\gam\in \Gam_\delta$, we have $\lim_nf(\gam(n))<\delta$,
while for almost every $\gam\in\Gam-\Gam_\delta$,
we have $\lim_nf(\gam(n))\geq \delta$.
Since both $\Gam_\delta$ and $\Gam-\Gam_\delta$
are non $p$-null, it follows that for every constant $c\in \R$
the set of $\gam\in \Gam$ with $\lim_nf(\gam(n))\neq c$
is not $p$-null.

Now theorem~\ref{hde} implies that there are non-constant,
$p$-harmonic functions,
on $G$.
\ep

\section{Sphere packing in $\R^d$ and a $d$-resolving metric}

In this section we will show that a quasi-sphere packing in $\R^d$ provides
a $d$-resolving metric for
$G$. The argument is similar to the case $d=2$ in ~\cite{bs:har}.
This will allow us to finish the proof of Theorem~\ref{packgen}.
\medskip

\begin{thm}
\label{pare}
Let $G$ be a graph with uniformly bounded vertex degree.
If $G$ admits a quasi-sphere packing in $\R^d$, then $G$ admits a $d$-resolving
metric.
\end{thm}

\noindent
{\bf Proof:} Assume there is a quasi-sphere packing of $G$ in $\R^d$, thus
it admits a
packing in $S^d$.
Take any $e\in E$, and let its vertices be $u,v$.
We set $m(e)=\diam(P_u)+\diam(P_v)$.
This defines a metric $m:E\to (0,\infty)$,
Because the packing $P$ is contained
in $S^d$, its total volume is finite,
and this implies that $m\in L^d(E)$.

We shall now show that $m$ is $d$-resolving.
For any $v\in V$ we let $z(v)$ denote the center of the
disk $P_v$. Let $p$ be any point in $\partial_m G$.
Let $v_1,v_2,\dots$ be a sequence in $V$ that converges to
$p$ in the completion of $G$, with respect to $m$, $C_m(G)$.
Then $\lim _{n,k\to\infty} d_m(v_n,v_k)=0$.
This easily implies that $\lim_{n,k\to\infty} |z(v_n)-z(v_k)|=0$.
We therefore conclude that the limit $\lim_n z(v_n)$ exists,
and denote this limit by $z(p)$.  If $w_1,w_2,\dots$ is another
sequence in $V$ that converges to $p$, then the limit
$\lim_n z(w_n)$ will still be $z(p)$.  This follows from
the fact that any ordering of the union $\{v_n\}\cup\{w_n\}$
as a sequence will still converge to $p$.  Hence $z(p)$ does not depend
on the sequence chosen.

Let $\Gam_p$ be the collection of all half-infinite paths
$\gam=(\gam(0),\gam(1),\dots)$ in $G$ that converge to $p$ in $C_m(G)$.
We need to show that $\Gam_p$ is $d$-null.
This will be done by producing an $L^d(E)$ metric $m_p$ such
that $\length_{m_p}(\gam)=\infty$ for every $\gam\in\Gam_p$.
The argument will be similar to an argument in \cite{sch:tel} and \cite{he-sc}.

We now inductively construct a sequence of positive numbers
$r_1>r_2>r_3>\dots$.
For $r>0$, let $B(r)$ denote the disk $\{z\in S^d:|z-z(p)|<r\}$.
Take $r_1=1$.  Suppose that $r_1,r_2,\dots ,r_{n-1}$ have
been chosen.  Let $r_n$ be in the range $0<r_n< r_{n-1}/2$
and be sufficiently small so that the two sets of vertices
$\{v\in V: z(v)\in B(2r_n)\}$ and $\{v\in V: z(v)\notin B(r_{n-1})\}$
are disjoint and there is no edge in $G$ connecting them.
To see that this can be done, observe that for any $r>0$
there are finitely many vertices $v\in V$ such that $\diam(P_v)\geq r$.
Since there is a uniform bound on the degrees, for every $r>0$ there is
a $\rho(r)\in(0,r/2]$ such that the closure of $B(\rho(r))$
does not intersect any of the sets $P_v$ satisfying $\diam(P_v)\geq r/2$.
This implies that there will be no $P_v$ that intersects both
circles $\partial B(r)$ and $\partial B(\rho(r))$.
Hence we may take $r_n = \rho(\rho(r_{n-1}))/2$.

For $r\in (0,\infty)$, let $\psi_r:\R^d\to [0,r]$ be defined by

\[ \psi_r(z)= \left\{ \begin{array}{ll}
                                    r     & \mbox{if } |z-z(p)|\leq r,\\
                               2r-|z-z(p)| & \mbox{if }r\leq |z-z(p)|\leq 2r,\\
                                   0 &\mbox{if }|z-z(p)|\geq 2r.
                                   \end{array}
                     \right. \]

In other words, $\psi_r$ is equal to $r$ on $B(r)$, equal to
$0$ outside $B(2r)$, and in the annulus $B(2r)-B(r)$ it is
linear in the distance from its center $z(p)$.
For each $n=1,2,\dots$ and $v\in V$, we define
$$
\phi_n(v) = \psi_{r_n}(z(v)).
$$
The construction of the sequence $r_1,r_2,\dots$ insures
that the supports of $d\phi_n$ and $d\phi_{n'}$ are disjoint when
$n\neq n'$.  It is easy to see that
the definition of $\phi_n$ shows that there is
a finite constant $C$ such that
\begin{equation}
\label{vf}
|d\phi_n(e)|^d\leq C\mbox{vol}((P_u\cup P_v)\cap B(3 r_n)),
\end{equation}
where $\{u,v\}= e$.
Let $\MD$ be an upper bound on the degrees of the vertices in
$G$.  Since the interiors of the sets in $P$ are disjoint, $(1)$
implies
$$
\|d\phi_n\|^d_d=\sum_{e\in E} d\phi_n(e)^d \leq 3^d \omega_d \pi C \MD r_n^d,
$$ %
where $\omega_d r^d$ is an upper bound for the volume of an $r$-ball in $S^d$. %

Now set
$$
\phi= \sum_{n=1}^\infty \frac{\phi_n}{n r_n}.
$$
As we have noted, the supports of the different $d\phi_n$ are
disjoint, and therefore,
$$
\|d\phi\|^d_d =\sum_{n=1}^\infty \frac{\|d\phi_n\|^d}{n^d r_n^d},
$$
and the above estimate for $\|d\phi_n\|^d_d$ shows that $|d\phi|\in L^d(E)$.
Therefore, there is some metric $m_p\in L^d(E)$ such that
$m_p (e)\geq |\phi_n(e)|$ for every $e\in E$.
(Technically, we cannot take $m_p=|d\phi|$,
since $|d\phi|$ is not positive, and hence not a metric.)

Now let $\gam=(\gam(1),\gam(2),\dots)$ be any path in $\Gam_p$,
and let $E(\gam)$ denote its edges.
We have $\lim_n z(\gam(n))=z(p)$.  Therefore,
$$
\lim_n\phi(\gam(n))=\lim_{z\to z(p)}\sum_j\frac{\psi_{r_j}(z)}{j r_j}
= \sum_j\frac 1j=\infty.
$$
This gives
$$
\length_{m_p}(\gam)=\sum_{e\in E(\gam)} m_p(e)
\geq \sum_{e\in E(\gam)}|d\phi(e)| =\infty.
$$
So $\Gam_p$ is $d$-null, and $m$ is $d$-resolving.
\ep
\medskip

\noindent
{\bf Proof of Theorem~\ref{packgen}:} Follows from Theorems \ref{resolving}
and \ref{pare}. %
\ep
\medskip

\noindent
{\bf Proof of Theorem~\ref{packZH}:} By Theorem~\ref{packgen} it is enough to
note that $\Z^{d+1}$, co-compact lattices in $\H^{d+1}$ and $3$-regular tree $\times \mbox{  } \Z$ are
$d$-nonparabolic, yet admit %
no non constant $d$-harmonic functions.

For
$\Z^{d+1}$, the first part can be found in \cite{ma} (see also \cite{ho:1}),
the second in \cite{ho-so:1}.

Hyperbolic space $\H^d$ has strictly positive Cheeger constant. This
property passes to co-compact lattices, see chapter 6 in \cite{glp}. Let
$G=(V,E)$ be the Cayley graph of a co-compact lattice in $\H^d$. By Cheeger's inequality,
$$\sum_{v\in V} |f(v)| \le c \sum_{e\in E} |df(e)|,$$
for $f$ with compact support.
Applying this to $f=u^p$, $p\ge 1$, yields
$$\sum_{v\in V} |f(v)|^p \le c \sum_{e\in E} |df(e)|^p .$$
This, in turn, implies that $G$
is not $p$-parabolic, see also~\cite{ho:3}. The same argument applies to $3$-regular tree $\times  \mbox{  } \Z$.

To deal with $p$-Dirichlet $p$-harmonic functions, we need some preliminaries.

If $f\in\ell^p (V)$ and $g\in\ell^{p'}(V)$, $p'=\frac{p}{p-1}$, let us denote by
\begin{eqnarray*}
\langle f,g\rangle =\sum_{v\in V} f(v)g(v).
\end{eqnarray*}
The $p$-Laplacian is related to the gradient of $D_p$ in the following way. If $\lap_p f\in \ell^{p'}(V)$ and $g\in\ell^p (V)$,
\begin{eqnarray*}
\frac{\partial}{\partial t}D_{p}(f+tg)_{t=0}=-p\langle g,\lap_p f\rangle.
\end{eqnarray*}
Since $D_p$ is homogeneous of degree $p$, differentiating $D_p (f+tf)=(1+t)^{p}D_{p}(f)$ leads to the following identity,
\begin{eqnarray*}
\langle f,\lap_p f\rangle=-D_p (f),
\end{eqnarray*}
for $f\in\ell^p (V)$ and $\lap_p f\in \ell^{p'}(V)$. In particular, a
non-constant $p$-Dirichlet $p$-harmonic function cannot be
$L^p$-integrable. Thus it suffices to show that for the graphs of interest,
every $p$-Dirichlet function is $L^p$-integrable, up to a constant. This
suggests the following definition. The degree $1$ exact $L^p$-cohomology (\emph{$L^p$-cohomology} for short) of a Riemannian manifold or a graph is the quotient of the space of $p$-Dirichlet functions by the subspace generated by constants and $L^p$-integrable $p$-Dirichlet functions.

We have just shown that constants are the only $p$-harmonic functions in the trivial cohomology class. More generally, if $p>1$, a cohomology class contains at most one $p$-harmonic function, up to constants. Indeed, since $D_p$ is strictly convex, if $g\not=0$, $g\in \ell^p (V)$ and $\lap_p f\in \ell^{p'}(V)$,
\begin{eqnarray*}
D_p (f+g)>D_p (f)-p\langle g,\lap_p f\rangle.
\end{eqnarray*}
So if $f$ is $p$-harmonic, $D_p (f+g)>D_p (f)$. If $f+g=h$ is $p$-harmonic as well, $D_p (f)=D_p (h-g)>D_p (h)$, contradiction.

The translation by one in the $\Z$ factor is an isometry $\tau$ of the graph $G=3$-regular tree $\times \mbox{  } \Z$. If $f$ is a function on vertices of $G$, and if $D_p (f)<\infty$, then $f\circ\tau-f$ is $L^p$-integrable, i.e. $f$ and $f\circ\tau$ belong to the same $L^p$-cohomology class. If $f$ is $p$-harmonic, so is $f\circ\tau$. This implies that $f\circ\tau=f$, i.e. $f$ is $\Z$-invariant. Since $D_p (f)<\infty$, this can hold only if $f$ is constant. So $p$-Dirichlet $p$-harmonic functions on $G$ are constant, for all $p>1$.

We use the rough isometry invariance of $L^p$-cohomology, not formally stated but proven in \cite{pa}, formally stated in chapter 8 of \cite{g}, see also \cite{bp}: if a bounded degree graph $G$ is roughly isometric to a bounded geometry Riemannian manifold $M$, then, for all $p\geq 1$, the $L^p$-cohomologies of $G$ and $M$ are isomorphic. A result of \cite{pa} states that $L^p$-cohomology of $\H^{d+1}$ vanishes if and only if $p\leq d$. It follows that co-compact lattices in hyperbolic space $\H^{d+1}$ have no non constant $d$-Dirichlet $d$-harmonic functions.
\ep

\medskip

\noindent
{\bf Remark:} In the course of the proof, we have shown that if a graph
admits non constant $p$-Dirichlet $p$-harmonic functions, then its
$L^p$-cohomology does not vanish. This yields the following corollary of Theorem \ref{packgen}.

\begin{cor}
\label{corpackgen}
If a bounded degree, $d$-nonparabolic graph can be quasi-sphere packed into
$\R^d$, then its $L^d$-cohomology does not vanish.
\end{cor}

Here is a consequence of Corollary \ref{corpackgen} in the world of
Cayley graphs of co-compact lattices in connected Lie
groups.

In \cite{pa-07}, it is shown that if a connected Lie
group has non vanishing $L^p$-cohomology for some $p$, then it is roughly isometric
either to a negatively curved Lie group,
or to a compact by solvable unimodular Lie group.
Thus a co-compact lattice in a connected Lie group which has non vanishing $L^p$-cohomology for some $p$ is either Gromov hyperbolic or amenable. In particular, lattices in higher rank semi-simple Lie groups cannot be quasi-sphere packed into $\R^d$ for any $d$.

\section{Further problems}

\noindent
{\bf Remark:} Holopainen and Soardi~\cite{ho-so:1} proved the Liouville
theorem for
$p$-harmonic functions on any $G$, for which a weak Poincar\'{e} %
inequality holds and has the volume doubling property.
This can be used %
to show that the Cayley graph of any group of %
polynomial volume growth with power $\geq d+1$ , cannot be (quasi-) sphere
packed in $\R^d$.

\begin{prob}
Can Cayley graphs of polynomial %
volume growth groups be %
packed in $\R^d$, for some $d$?
\end{prob}

\begin{prob}
Show that any sphere packing of $\Z^3$ in $S^3$, has only %
one accumulation point
for the centers of the spheres.
\end{prob}

\begin{prob}
By compactness it follows, that already a sufficiently large ball in
$\Z^{d+1}$, cannot be packed in $\R^d$. Estimate the size of such a ball. %
\end{prob}

\begin{prob}
Characterize or give necessary and sufficient conditions regarding which bounded geometry %
Riemannian metrics on $\R^3$
admit a rough isometric graph, packable in $\R^3$. %
\end{prob}

\begin{prob}
Packing in other spaces. Separation can be used to show that $\Z^d$ can not be packed in $\H^d$. What about packing of $\Z^3$
or co compact lattices of $\H^3$ in $\H^2 \times \R$? 
\end{prob}

The recent paper~\cite{bc} uses Theorem~\ref{pare} above in order to extend %
parts of ~\cite{bs:lim} to higher dimensions
and suggests further related problems as well.
\medskip

\noindent
{\bf Acknowledgements.} Thanks to I. Holopainen
and L. Saloff-Coste for a useful discussion. Special thanks to
P. Pansu for writing the last part using $L^p$-cohomology.


\noindent

\begin{thebibliography} {YMN}

\bibitem{bc}
Benjamini I. and Curien, N. On limits of graphs sphere packed in Euclidean space and applications.
{\it European J. Comb.} To appear.
arXiv:0907.2609

\bibitem{bs:har} Benjamini, I. and Schramm, O. Harmonic functions on planar
and almost planar graphs and manifolds, via circle packings.
{\it Invent. Math.} {\bf 126} (1996), 565--587

\bibitem{bs:lim}
Benjamini I. and Schramm, O. Recurrence of Distributional Limits of Finite Planar Graphs.
{\it Electron. J. Probab.} {\bf 6}  (2001) 13 pp.

\bibitem{bs:sep} Benjamini, I. Schramm, O. and Tim\'ar, A.  Separation.
(2009), {\it preprint}
 arXiv:1004.0921

\bibitem{bp} Bourdon, M. and Pajot, H. Cohomologie $l_p$ et espaces de Besov.
{\it J. Reine Angew. Math.} {\bf 558}  (2003), 85--108 %

\bibitem{b-s} Brightwell, G.R. and Scheinerman, E.R. Representation of
planar graphs. {\it SIAM J. Discrete Math.} {\bf 6} (1993), 214--229

\bibitem{co} Colin De Verdi\`{e}re, %
Y. Un principe variationnel pour les
empilements de cercles. {\it Invent. math.} {\bf 104} (1991), 655--669

\bibitem{co-sl} Conway, J. H. Sloane, N. J. A. Sphere packings, Lattices and
Groups. {\it Springer-Verlag}, New-York, (1988)


\bibitem{du} Duffin, R. J. The extremal length of a network.
{\it J. Math. Anal. Appl.} {\bf 5} (1962), 200--215 %

\bibitem{glp} Gromov, M. Metric structures for Riemannian and
  non-Riemannian spaces. Based on the 1981 French original. With appendices
  by M. Katz, P. Pansu and S. Semmes. Translated from the French by Sean
  Michael Bates. Progress in Mathematics, {\bf 152}. Birkh\"a user Boston, Inc., Boston, MA, 1999.

\bibitem{g} Gromov, M. Asymptotic invariants of infinite groups. Geometric group theory, Vol. 2 (Sussex, 1991), 1--295, London Math. Soc. Lecture Note Ser., {\bf 182}, Cambridge Univ. Press, Cambridge, (1993)

\bibitem{he-sc} He, Z.-X. and Schramm, O. Hyperbolic and parabolic packings.
{\it Discrete Comput. Geom.} {\bf 14} (1995), 123--149

\bibitem{ho:1} Holopainen, I. Positive solutions of quasilinear elliptic
equations on Riemannian manifolds. {\it Proc. London Math. Soc.} {\bf 65}
(1992), 651--672

\bibitem{ho:2} Holopainen, I. Rough isometries and $p$-harmonic functions
with finite Dirichlet integral. {\it Rev. Mat. Iberoamericana}
{\bf 10} (1994), 143--176

\bibitem{ho:3} Holopainen, I. Volume growth, Green's functions, and
parabolicity of ends. (1997), {\it preprint}

\bibitem{ho-so:1} Holopainen, I. and Soardi, M. P. A strong Liouville theorem
for $p$-harmonic functions on graphs. {\it Ann. Acad. Sci. Fennicae}
{\bf 22} (1997), 205--226

\bibitem{ho-so:2} Holopainen, I. and Soardi, M. P. $p$-harmonic functions on
graphs and manifolds. {\it Manuscripta Math.} {\bf 94} (1997), no. 1, 95--110 %

\bibitem{ko} Koebe, P. Kontaktprobleme der konformen Abbildung.
{\it Berichte Verhande. S\"{a}chs. Akad. Wiss. Leipzig, Math.-Phys. Klasse}
{\bf 88} (1936), 141--164

\bibitem{ku-sc} Kuperberg, G. and Schramm, O. Average kissing numbers for
non-congruent sphere packings. {\it Math. Res. Lett.} {\bf 1}
(1994), 339--344

\bibitem{ma} Maeda, F.Y. A remark on the parabolic index of infinite networks.
{\it Hiroshima Math. J.} {\bf 7} (1977), 147--152

\bibitem{m-r} Marden, A. and Rodin, B. On Thurston's formulation and proof of
Andreev's theorem. In Computational Methods and Function Theory,
{\it LNM 1435}, (1989), 103--115

\bibitem{pa} Pansu, P. Cohomologie $L^p$ des vari\'et\'es \`a courbure n\'egative,
cas du degr\'e $1$.  %
{\it Rend. Sem. Mat. Univ. Politec. Torino 1989}, Special Issue, (1990),
95--120.

\bibitem{pa-07} Pansu, P.
Cohomologie $L^{p}$ en degr\'e $1$ des espaces homog\`enes.
{\it Potential Anal.} {\bf 27}, (2008), 151--165

\bibitem{sch:tel} Schramm, O. Transboundary extremal length. {\it J.
d' Analyse Math.} {\bf 66} (1995), 307--329

\bibitem{so:1} Soardi, M. P. Rough isometries and Dirichlet finite harmonic
functions on graphs. {\it Proc. AMS} {\bf 119} (1993), 1239--1248

\bibitem{so:2} Soardi, M. P. Potential theory on infinite networks.
{\it LNM 1590} Springer, (1994)

\bibitem{so-ya} Soardi, M. P.  and Yamasaki, M. Parabolic index and rough
isometries. {\it Hiroshima Math. J.} {\bf 23} (1993), 333--342

\bibitem{ya:1} Yamasaki, M. Parabolic and hyperbolic infinite networks.
{\it Hiroshima Math. J.} {\bf 7} (1977), 135--146

\bibitem{ya:2} Yamasaki, M. Ideal boundary limit of Dirichlet functions.
{\it Hiroshima. Math. J.} {\bf 16} (1986), 353--360

\end{thebibliography}
\end{document}